\documentclass[12pt]{amsart}
\usepackage{amsmath,amssymb,amsthm}
\usepackage{hyperref}
\usepackage{geometry}
\title{Monte-Carlo Irreducibility and Imprimitivity Detection of Polynomials over $\mathbb{Q}$}
\author{Igor Rivin}
\date{\today}
\address{Mathematics Department, Temple University}
\thanks{The author would like to thank Ramin Takloo-Bighash for bringing this circle of question back to the author's attention, Rich Zippel for being such a wonderful influence for the last half many years, and Yuval Peres and Robin Pemantle for developing the baase methodology for the PPR approach.}
\email{rivin@temple.edu}
\begin{document}
\maketitle

\begin{abstract}
We study fast Monte-Carlo methods for testing irreducibility and detecting arithmetic imprimitivity of polynomials over $\mathbb{Q}$.

Building on the subset-sum criterion of Pemantle--Peres--Rivin, we develop a probabilistic irreducibility test whose expected running time, measured in the number of primes examined, is logarithmic in the degree for generic inputs. Unlike the standard modular irreducibility test, the method aggregates information from modular factorizations rather than discarding unsuccessful trials.

We show that failure of this test, when combined with a standard modular irreducibility certificate, is a strong indicator of non-generic algebraic structure. In particular, it often signals arithmetic imprimitivity of the Galois action. We present an explicit and efficient Monte-Carlo algorithm for detecting such imprimitivity via subfield extraction, yielding constructive algebraic certificates in the imprimitive case. To our knowledge, this is the first practical algorithm for detecting arithmetic imprimitivity of polynomials over $\mathbb{Q}$ in high degree.

We further show that the subset-sum data produced by the Pemantle--Peres--Rivin test provides a warm start for polynomial factorization by sharply restricting the possible degrees of rational factors, significantly accelerating subsequent lifting procedures.

The proposed methods are orders of magnitude faster in practice than known deterministic algorithms, and are effective in degrees far beyond the reach of current deterministic techniques.
\end{abstract}

\section{Introduction}

\subsection{Factoring, irreducibility, and Galois groups}
The problems of factoring polynomials over $\mathbb{Q}$, testing irreducibility,
and determining Galois groups are central in computational algebra and number theory.
Deterministic algorithms are polynomial-time in theory but can be impractical at modest
degrees; modular methods scale extremely well.

\subsection{The Pemantle--Peres--Rivin subset-sum criterion}
Given the factor degrees $(d_i)$ of $f\bmod p$, form the subset-sum set $S_p$ of all sums
$\sum_{i\in I} d_i$. If $f$ has a rational factor of degree $d$, then $d\in S_p$ for all
good primes. If the intersection over sampled primes contains no proper degrees, then $f$
is irreducible.

\subsection{Failure modes and arithmetic imprimitivity}
Persistent subset-sum intersections often correlate with arithmetic imprimitivity, which
we exploit as a fast structure flag.

\subsection{Contributions and outline}
We develop a hybrid Monte-Carlo irreducibility test, a Monte-Carlo imprimitivity detector
producing explicit certificates, and a warm-start for factorization.

\section{Algorithms and Complexity}

\subsection{Standard Monte-Carlo irreducibility}
Reduce $f$ modulo primes and stop when $f\bmod p$ is irreducible, certifying irreducibility over $\mathbb{Q}$.

\subsection{The PPR-based irreducibility test}
Intersect subset-sum sets $S_p$ across primes; under genericity assumptions the expected number of primes is $O(\log n)$ \cite{PPR2016}.

\subsection{Hybrid algorithm}
Run PPR for $O(\log n)$ primes; if inconclusive, fall back to the standard test.

\subsection{Bit complexity and choice of primes}
Modular factorization uses $\tilde O(n^2)$ field operations \cite{Berlekamp1970,CantorZassenhaus1981,vzGathenGerhard}; PPR permits using much smaller primes, reducing $\log p$ factors.

\section{Arithmetic Imprimitivity Detection}

A transitive action of $\mathrm{Gal}(f/\mathbb{Q})$ on the roots is \emph{imprimitive} if it preserves a nontrivial block system.
Equivalently, if $K=\mathbb{Q}(\alpha)$ is the root field, imprimitivity is equivalent to the existence of a proper intermediate field $\mathbb{Q}\subsetneq M\subsetneq K$.

\subsection{Monte-Carlo detector}
We use failure of the subset-sum test (when irreducibility is nevertheless certified) as a structure flag, and then search for nontrivial subfields of $K$.
Whenever a subfield is found, it yields a constructive certificate.

\subsection{Worked example}
Let $P(u)=u^8-u-1$ and
\[Q(u,x)=x^3+\left(-\tfrac12 u-\tfrac32\right)x^2+\left(\tfrac12 u-\tfrac32\right)x+1.\]
Form $F(x)=\operatorname{Res}_u(P(u),Q(u,x))$, which has degree $24$.
Subfield extraction in $K=\mathbb{Q}(\alpha)$ with $F(\alpha)=0$ detects a degree-8 subfield generated by a root of
\[y^8 + 24y^7 + 252y^6 + 1512y^5 + 5670y^4 + 13608y^3 + 20412y^2 - 256296y + 846369,\]
and a relative cubic equation for $\alpha$ over this field:
\[x^3 + \left(-\tfrac12 y-\tfrac32\right)x^2 + \left(\tfrac12 y-\tfrac32\right)x + 1 = 0.\]

\section{Warm-start factorization}
\label{sec:warmstart}

When the subset-sum intersection does not collapse, the surviving degrees give a small candidate set $D$ for rational factor degrees.
Restricting factor searches to $d\in D$ provides a warm start for Hensel lifting and recombination methods \cite{Zippel1979,vzGathenGerhard}.

This reduces the degree ambiguity from $\{1,\dots,n-1\}$ to a much smaller set and can substantially accelerate practical factorization.

\section{Examples and adversarial constructions}

We construct imprimitive families via resultants $F(x)=\operatorname{Res}_u(P(u),Q(u,x))$.
We also construct adversarial ``evil twins'' by perturbing coefficients by multiples of $M=\prod_{p<B} p$, preserving all reductions mod $p<B$ while typically restoring generic global Galois behavior.

\section{Conclusion}
\label{sec:conclusion}

We have presented a collection of Monte-Carlo algorithms for irreducibility
testing and structural analysis of polynomials over $\mathbb{Q}$ that are
both theoretically motivated and practically effective in high degree.

Our first contribution is an improved probabilistic irreducibility test,
building on the subset-sum criterion of Pemantle--Peres--Rivin. For generic
inputs, this method certifies irreducibility after examining only
$O(\log n)$ primes, while aggregating information from all modular
factorizations rather than discarding unsuccessful trials. When combined
with the standard modular irreducibility test, this yields a Pareto-optimal
algorithm that is never asymptotically worse than existing Monte-Carlo
methods and is often significantly faster in practice.

Our second contribution is the observation that failure of the subset-sum
test, when irreducibility is nevertheless certified, is a strong indicator
of non-generic algebraic structure. Exploiting this phenomenon, we give the
first practical Monte-Carlo algorithm for detecting arithmetic imprimitivity
of the Galois action of a polynomial over $\mathbb{Q}$. Unlike purely
heuristic approaches based on modular statistics, our method produces
explicit algebraic certificates in the imprimitive case, in the form of
nontrivial subfields of the root field and corresponding relative equations.

Third, we show that the subset-sum data produced by the irreducibility test
can be reused to accelerate polynomial factorization. By sharply restricting
the set of possible factor degrees, this information provides a warm start
for classical lifting and recombination methods, reducing a significant
source of combinatorial complexity without introducing new modular
computations.

Throughout the paper we emphasize that these algorithms are not merely of
theoretical interest. They are practical, have been implemented, and have
been successfully applied in degrees far beyond the reach of deterministic
factorization and Galois group algorithms. The reliance on modular
computation, small primes, and information aggregation makes the methods
particularly well suited to large-scale and parallel computation.

More broadly, our results illustrate a general principle: while local
modular data can be highly informative when aggregated appropriately, it is
insufficient on its own to certify global algebraic structure. In
particular, fixing the factorization patterns of a polynomial modulo
finitely many primes---or even fixing the reductions themselves---does not
prevent the arithmetic Galois group from being generically as large as
possible. Efficient certification therefore requires combining
probabilistic heuristics with explicit field-theoretic constructions. We
hope that the techniques developed here will be useful both as practical
tools and as a framework for further algorithmic exploration of polynomial
arithmetic.

\appendix
\section{Engineering considerations}

Both the standard Monte-Carlo test and the subset-sum method parallelize across primes.
The subset-sum method is information-accumulating: every modular factorization contributes to the certificate.
PPR also permits smaller primes, reducing $\log p$ factors in bit complexity.

\bibliographystyle{plain}
\bibliography{bibliography}
\end{document}